\documentclass{amsart}

\renewcommand{\oddsidemargin}{0cm}

\newtheorem{theorem}{Theorem}[section]
\newtheorem{lemma}[theorem]{Lemma}

\newtheorem{proposition}[theorem]{Proposition}

\theoremstyle{definition}
\newtheorem{definition}[theorem]{Definition}
\newtheorem{example}[theorem]{Example}
\newtheorem{remark}[theorem]{Remark}

\usepackage{amssymb}

\def\R{\mathbb R}
\def\X{\mathbb X}
\def\G{\mathbb G}
\def\J{\mathbb J}

\def\cal{\mathcal}

\def\Lip{\ensuremath{\textsl{Lip}}}
\def\lip{lip}

\def\ds{\displaystyle}
\def\R{\mathbb R}

\def\cal{\mathcal}

\newcommand{\defnemph}[1]{\textbf{#1}}
\newcommand{\Xts}{\ensuremath{\{X(t,s)\}_{(t,s)\in\Delta}} }
\newcommand{\Uts}{\ensuremath{\{U(t,s)\}_{(t,s)\in\Delta}} }

\begin{document}

\title[On the asymptotic behavior of the solutions of semilinear nonautonomous equations]
%{On the admissibility of function spaces and the exponential stability of nonlinear evolutionary processes}
%{Stability of solutions of semilinear evolution equations using the ``test functions method''}
{On the asymptotic behavior of the solutions of semilinear nonautonomous equations}

\author{Nguyen Van Minh}
\address{Department of Mathematics and Philosophy, Columbus State University, Columbus, GA 31907, USA}
\email{nguyen\_minh2@columusstate.edu}

\author{Gaston M. N'gu\'er\'ekata}
\address{Department of Mathematics, Morgan State University, Baltimore, Maryland 21251, USA}
\email{Gaston.N'Guerekata@morgan.edu}

\author{Ciprian Preda}
\address{Department of Mathematics, West University of Timisoara, Romania}
\email{ciprian.preda@feaa.uvt.ro}

\begin{abstract}

We consider nonautonomous semilinear evolution equations of the form
\begin{equation}
\label{semilineq}
\frac{dx}{dt}= A(t)x+f(t,x) \ .
\end{equation}
Here $A(t)$ is a (possibly unbounded) linear operator acting on a real or
complex Banach space $\X$ and $f: \R\times\X\to\X$ is a (possibly nonlinear)
continuous function. We assume that the linear equation \eqref{lineq} is well-posed (i.e. there exists a
continuous linear evolution family \Uts such that for every $s\in\R_+$ and $x\in D(A(s))$,
the function $x(t) = U(t, s) x$ is the uniquely determined solution
of equation \eqref{lineq} satisfying $x(s) = x$). Then we can consider the \defnemph{mild solution}
of the semilinear equation \eqref{semilineq} (defined on some interval $[s , s + \delta), \delta > 0$)
as being the solution of the integral equation
\begin{equation}
\label{integreq}
x(t) = U(t, s)x + \int_s^t U(t, \tau)f(\tau, x(\tau)) d\tau \quad,\quad t\geq s\ ,
\end{equation}
 Furthermore, if we assume also that the nonlinear function $f(t, x)$ is jointly continuous with respect to $t$ and $x$ and
Lipschitz continuous with respect to $x$ (uniformly in $t\in\R_+$, and $f(t,0) = 0$ for all $t\in\R_+$) we can generate a (nonlinear) evolution family \Xts, in the sense that the map $t\mapsto X(t,s)x:[s,\infty)\to\X$
is the unique solution of equation \eqref{integreq}, for every $x\in\X$ and $s\in\R_+$.

 Considering the \emph{Green's operator} $(\G f)(t)=\int_0^t X(t,s)f(s)ds$ we prove that if the following conditions hold
\begin{itemize}
        \item \quad the map $\G f$ lies in $L^q(\R_+,\X)$ for all $f\in L^{p}(\R_+,\X)$, and
        \item \quad $\G:L^{p}(\R_+,\X)\to L^{q}(\R_+,\X)$ is Lipschitz continuous, i.e. there exists $K>0$ such that
$$\|\G f-\G g\|_{q} \leq K\|f-g\|_{p}\ ,\ \mbox{for all}\ f,g\in L^p(\R_+,\X)\ ,$$
\end{itemize}
then the above mild solution will have an exponential decay.

\end{abstract}

\subjclass[2010]{Primary 34D05,34G10; Secondary 47D06, 93D20}
\keywords{semilinear evolution equations,   exponential stability}

\date{\today}

\maketitle

\section{Introduction}

In the last few decades, we can note an increasing research interest in the asymptotic behavior
of the solutions of the linear differential equations
\begin{gather}
\label{lineq}
\frac{dx}{dt}=A(t)x\ ,\ t\in[0,\infty)\ ,\ x\in\X\ ,
\end{gather}
where $A(t)$ is in general an unbounded linear operator on a Banach space $\X$,
for every fixed $t$.

In the case that $A(t)$ is a matrix continuous function, O.~Perron \cite{per} first observed
a connection between the asymptotic behavior of the solutions of the above equation and the properties
of the differential operator $\frac{d}{dt} - A(t)$ as an operator on the space
$C_b(\R_+,\R^n)$ of $R^n$-valued, bounded and continuous functions on the half-line $\R_+$.

This result became a milestone for many works on the
qualitative theory of solutions of differential equations.
We refer the reader to the monograph by Massera and Sch\"{a}ffer \cite{masseraschaffer},
and Daleckij and Krein \cite{dalkre} for a characterization of the exponential
dichotomy of the solutions of the above equation in terms of surjectiveness
of the differential operator $\frac{d}{dt} - A(t)$ in the case of bounded $A(t)$ and by
Levitan and Zhikov \cite{levzhi} for an extension to the infinite-dimensional setting
for equations defined on the whole line.

Recently, there is a lot of work done in the study of the asymptotic behavior
of the solutions of differential equations in Banach spaces, in particular in
the unbounded case (see e.g.,
\cite{EngelNagel,Neervenbook,prepogpre3,prepogpre5,schnaubelt2001}).

Also, there is an increasing interest ``applications-wise''
for the Cauchy problems associated with the above evolution equations
since many physical situations can be interpreted as Cauchy
problems by choosing the ``right'' state space. Among the physical
equations which are eligible for such an approach we can mention
heat equation, Schr\"odinger equation, certain population
equations, Maxwell's equations, wave equation, and also
seemingly unrelated problems such as delay equations, Markov
processes or Boltzmann's equations.

In the present paper, we will continue the approach initiated by Perron
for semilinear nonautonomous evolution equations of the form
\begin{equation}
\label{semilineq}
\frac{dx}{dt}= A(t)x+f(t,x) \ .
\end{equation}
Here $A(t)$ is a (possibly unbounded) linear operator acting on a real or
complex Banach space $\X$ and $f: \R\times\X\to\X$ is a (possibly nonlinear)
continuous function. Following \cite{aulmin}, we assume that the linear equation \eqref{lineq} is well-posed (i.e. there exists a
continuous linear evolution family \Uts such that for every $s\in\R_+$ and $x\in D(A(s))$,
the function $x(t) = U(t, s) x$ is the uniquely determined solution
of equation \eqref{lineq} satisfying $x(s) = x$). Then we can consider the \defnemph{mild solution}
of the semilinear equation \eqref{semilineq} (defined on some interval $[s , s + \delta), \delta > 0$)
as being the solution of the integral equation
\begin{equation}
\label{integreq}
x(t) = U(t, s)x + \int_s^t U(t, \tau)f(\tau, x(\tau)) d\tau \quad,\quad t\geq s\ ,
\end{equation}
 Furthermore, if we assume also that the nonlinear function $f(t, x)$ is jointly continuous with respect to $t$ and $x$ and
Lipschitz continuous with respect to $x$ (uniformly in $t\in\R_+$, and $f(t,0) = 0$ for all $t\in\R_+$) we can generate a (nonlinear) evolution family \Xts, in the sense that the map $t\mapsto X(t,s)x:[s,\infty)\to\X$
is the unique solution of equation \eqref{integreq}, for every $x\in\X$ and $s\in\R_+$.

 Considering the \emph{Green's operator} $(\G f)(t)=\int_0^t X(t,s)f(s)ds$ we prove that if the following conditions hold
\begin{itemize}
        \item \quad the map $\G f$ lies in $L^q(\R_+,\X)$ for all $f\in L^{p}(\R_+,\X)$, and
        \item \quad $\G:L^{p}(\R_+,\X)\to L^{q}(\R_+,\X)$ is Lipschitz continuous, i.e. there exists $K>0$ such that
$$\|\G f-\G g\|_{q} \leq K\|f-g\|_{p}\ ,\ \mbox{for all}\ f,g\in L^p(\R_+,\X)\ ,$$
\end{itemize}
then the above mild solution will have an exponential decay.

It is worth to note that, although the autonomous case
(i.e. time invariant evolution equations),
was  much more analyzed than the  nonautonomous case, the latter one
often arises quite naturally, not only in physics and mechanics,
but also in the mathematical theory of differential equations when one linearizes an autonomous equation along
a nonstationary solution.
For particular cases of autonomous evolution
equations arising from the linearization along a compact invariant manifold it
has been shown (see e.g. \cite{sacsel}) that one can define a skew-product semiflow
which allows to apply the methods of classical dynamical systems to the
underlying time-dependent equations.

For the case of a time-invariant linear part of equation \eqref{semilineq} the existence
problem for solutions has been investigated by many authors
(see e.g. \cite{iwamiya,komatsu,martin,oharu,oharu2,pavel,pazy,webb}  and the references therein).

\section{Semilinear evolution equations. Examples}
First, let us recall some notations and definitions. Throughout
this paper, $\X$ will denote a Banach space, $\R$ the set of all real numbers,
$\R_+$ the subset of all nonnegative real numbers and put $\Delta=\{(t,s)\in\R^2_+ : t\geq s\}$.
If $\mathbb{Y}$ denotes also a Banach space, then the set of all maps $T:\X\to\mathbb{Y}$ such that
$$\|T\| _{\lip} := \inf\{ M > 0 : \| Tx-Ty\| \leq M\| x-y \|\ ,\ \mbox{for all}\ x, y\in \X\}\ <\ \infty\ .$$
will be denoted by $\Lip(\X,\mathbb{Y})$.
Also, if $\X=\mathbb{Y}$ we will put simply $\Lip(\X)$ instead of $\Lip(\X,\X)$.
It is easy to see that $(\Lip(\X), \| \cdot \| _{\lip})$ is a
seminormed vector space which has the property
$$\| T \circ S \|_{\lip} \leq \| T \| _{\lip}\| S \| _{\lip}\ \mbox{for all} \ T,S \in\Lip(\X)\ .$$
For a given interval $\mathbb{J}$ of the real line, we denote by
$$L^{p}(\J,\X) = \{ f: \J\rightarrow \X: f\ \mbox{is measurable and}\ \int_{\J}\! \| f(t)\|^{p}dt < \infty \}\quad ,$$
for all $p \in [1,\infty)$ and by
$$
L^{\infty}(\J,\X)= \{ f:\J\rightarrow \X\,:\,f\ \mbox{is measurable and}\ {\rm ess}\sup\limits_{\hspace{-4mm}t\in \J}\| f(t)\| < \infty \} .
$$
It is well-known that $L^{p}(\J,\X), L^{\infty}(\J,\X)$ are Banach
spaces endowed with the norms
$$\| f\| _{p} = \left(\int_{\J}\| f(t)\|^{p}dt\right)^{1/p}\quad,$$
$$\| f\| _{\infty} = {\rm ess}\sup\limits_{\hspace{-4mm}t\in\J} \| f(t)\| \quad,$$
respectively.

\begin{definition}
\label{defn:evolproc}
A family \Xts of (possibly nonlinear) operators acting on $\X$
is called an \defnemph{evolution family} if it satisfies the following conditions:
\begin{enumerate}
\item[$(e_1)$] \quad $X(t,t)x = x$ for all $t\geq 0$ and $x \in \X;$
\item[$(e_2)$] \quad $X(t,s) = X(t,r)\circ X(r,s)$ for all $t\geq r \geq s \geq 0$.
\end{enumerate}
Such an evolution family is called \defnemph{continuous} if
there exist $M, \omega> 0$ such that
\begin{enumerate}
\item[$(e_3)$] \quad $\|X(t,s)\|_{\lip}\leq Me^{\omega(t-s)}$
\item[$(e_4)$] \quad $X(t,s)x$ is jointly continuous with respect to $t$, $s$ and $x$..
\end{enumerate}
\end{definition}

If the operators $X(t,s)$ are linear, then
by $(e_3)$ they are also bounded and the family will be called a \defnemph{linear evolution family}.

\begin{remark}
Condition $(e_3)$ is equivalent with the existence of some locally bounded function $\varphi$ such that
\begin{enumerate}
\item[$(e'_3)$] $\|X(t,s)x-X(t,s)y\|\leq \varphi(t-s)\|x-y\|$ for all $x,y\in \X$.
\end{enumerate}

Indeed, if $(e_3)$ holds one can take $\varphi(t)=Me^{\omega t}$. Conversely, the constants
$M=\sup_{t\in[0,1]}\varphi(t)$ and $\omega=\max\{1,\ln\varphi(1)\}$ satisfy $(e_3)$.
\end{remark}

\begin{definition}
\label{defn:wellposedness}
The linear equation \eqref{lineq} is said to be \defnemph{well-posed}
if there exists a continuous linear evolution family \Uts
such that for every $s\in\R_+$ and $x\in D(A(s))$,
the function $x(t) = U(t, s) x$ is the uniquely determined solution
of equation \eqref{lineq} satisfying $x(s) = x$.
\end{definition}

\begin{definition}
Suppose the linear equation \eqref{lineq} is well-posed.
Then, every solution $x(t)$ (defined on some interval $[s , s + \delta), \delta > 0$)
of the integral equation
\begin{equation}
\label{integreq}
x(t) = U(t, s)x + \int_s^t U(t, \tau)f(\tau, x(\tau)) d\tau \quad,\quad t\geq s\ ,
\end{equation}
is called a \defnemph{mild solution} of the semilinear equation \eqref{semilineq}
starting from $x$ at $t = s$.
Furthermore, \defnemph{equation \eqref{semilineq} is said to generate an evolution family} \Xts if
for every $x\in\X$ and $s\in\R_+$, the map $t\mapsto X(t,s)x:[s,\infty)\to\X$
is the unique solution of equation \eqref{integreq}.
\end{definition}

\begin{proposition}
Suppose the following conditions are satisfied:
\begin{enumerate}
\item[(i)] The linear equation \eqref{lineq} is well-posed;
\item[(ii)] The nonlinear function $f(t, x)$ is jointly continuous with respect to $t$ and $x$ and
        Lipschitz continuous with respect to $x$, uniformly in $t\in\R_+$, and $f(t,0) = 0$ for all $t\in\R_+$.
\end{enumerate}
Then, the semilinear equation \eqref{semilineq} generates a continuous evolution family.
\end{proposition}
\begin{proof}
Using standard arguments, see for instance \cite{segal} it can be shown that the equation
\eqref{semilineq} generates an evolution family. Moreover, from \cite{segal}, it follows that
$X(t,s)x$ is jointly continuous with respect to $t$, $s$ and $x$.
We show briefly bellow that $X(t,s)x$ also fulfill condition $(e_3)$ in Definition~\ref{defn:evolproc}.
We have that
\begin{gather*}
\begin{split}
\|X(t,s)x -X(t,s)y\|\leq&\, \|U(t,s)x-U(t,s)y\|+ \int_s^t \|U(t,\xi)\|\,\|f(\xi,X(\xi,s)x)-f(\xi,X(\xi,s)y)\|d\xi\\
  \leq&\, Ke^{\omega(t-s)}\|x-y\| + \int_s^t Ke^{\omega(t-\xi)}L\|X(\xi,s)x-X(\xi,s)y\|d\xi\ ,
\end{split}
\end{gather*}
where $L$ is a Lipschitz constant of $f(t,x)$ with respect to $x$ and $K,\omega$
stem from the well-posedness of the linear equation \eqref{lineq} and for convenience choose $\omega$ to be positive.
Applying Gronwall's Lemma we get
\begin{gather*}
\|X(t,s)x-X(t,s)y\|\leq Ke^{(\omega+KL)(t-s)}\|x-y\|\ ,
\end{gather*}
for any $(t,s)\in\Delta$ and $x,y\in X$.

\end{proof}

\begin{remark}
 It is worth to note that one of the goals with respect to the asymptotic behavior of solutions of the
equation \eqref{semilineq} is also to point out conditions for that equation to admit invariant
(stable, unstable or center) manifolds
(see, e.g., \cite{aulmin,batesjones,dalkre,halemagalhaes,henry,hps,minh,sellyou}).
As far as we know, the most popular conditions for the existence of invariant manifolds are the exponential stability, dichotomy
of the linear part \eqref{lineq}
and the uniform Lipschitz continuity of the nonlinear part $f(t, x)$
with sufficiently small Lipschitz constants (i.e., $\|f(t, x)-f(t, y)\|\leq M \|x-y\|$ for $M$ small enough).
Moreover, the manifolds considered in the existing literature
are mostly constituted by trajectories of solutions bounded on the positive (or negative)
half-line. We refer the reader to \cite{aulmin,batesjones,halemagalhaes,henry,hps,minh,sellyou}
and references therein for more details on this topic.

\end{remark}

\begin{example}
\label{firstex}
For a given  continuous map $h:{{\R } }\rightarrow[1/2,1]$ consider the differential equation
on $\R$
\begin{equation}
\label{(2.1.)}
 \dot{u}(t) = h(u(t))\ .
\end{equation}
We claim that this equation leads to a continuous evolution family.

Indeed, consider the map $H: {{\R } }\rightarrow {{\R } }$ given by
$$H(t)=\int^{t}_{0}\frac{ds}{h(s)}\quad .$$
One can easily check that
$$\vert u-v\vert \leq  \vert H(u)-H(v)\vert \leq 2\vert u-v\vert$$
for all $u,v\in\R$.
It follows that $H$ is bijective and so it is easy to see that
\begin{equation}
X(t,s):\R\to\R \quad,\quad X(t,s)x = H^{-1}(t-s+H(x))
\end{equation}
is an evolution family which has the property $(e_4)$.
Also, we have
\begin{gather*}
\begin{split}
\vert X(t,s)x-X(t,s)y\vert =& \vert H^{-1}(t-s+H(x))-H^{-1}(t-s+H(y))\vert \\
    \leq & \vert (t-s+H(x)) - (t-s+H(y))\vert \\
    \leq & \vert x-y\vert \quad ,
\end{split}
\end{gather*}
for all $(t,s)\in\Delta$ and all $x\in\R$.
Thus we obtain that  \Xts is a continuous evolution family on $\R$.
\end{example}

\begin{example}
\label{secondex}
Consider $h:\R\to\R$ continuously differentiable with $h^{\prime}\in L^{\infty}(\R,\R)$ and the problem
\begin{gather}
\label{(2.5.)}
\begin{cases}
\ds\frac{\partial u}{\partial t}(x,t)=\frac{\partial ^{2}u}{\partial x^{2}}(x,t)+h(u(x,t))\\[3mm]
\ds\frac{\partial u}{\partial x}(0,t)=\frac{\partial u}{\partial x}(1,t)=0
\end{cases}
\end{gather}
If we denote $x(t)=u(\,\cdot\,,t)$,  the problem \eqref{(2.5.)} is equivalent to
\begin{equation}\label{(2.6.)}
  \dot{x}(t)=Ax(t)+f(x(t)),
\end{equation}
where $A:D(A)\subset L^{2}([0,1],\R)\to L^{2}([0,1],\R)$,
$D(A)$ is defined as the set of all functions $z\in L^{2}([0,1],\R)$ such that
$z, z^{\prime}$ are absolutely continuous with $z^{\prime \prime}\in L^{2}([0,1],\R)$
and $z^{\prime}(0)=z^{\prime}(1) = 1$, and
for each $z\in D(A)$, we define $Az=z^{\prime \prime}$.

Also we consider $B:L^{2}([0,1],{{\R } })\rightarrow
L^{2}([0,1],{{\R } })$ which is given by $Bz=h\circ z$. It is well
known that $A$ generates a strongly continuous semigroup of linear
operators $\{ T(t)\}_{t\geq 0}$ in $L^{2}([0,1],{{\R } })$. And it
is easy to check that $B$ is Lipschitz continuous.

By Example~\ref{firstex} it is clear
that the equation (\ref{(2.6.)}), as an abstract variant of
(\ref{(2.5.)}), generates a  continuous evolution family.
\end{example}

For \Xts an evolution family, the \emph{trajectory} determined by $t_0\in\R_+$ and $x_0\in\X$ will be denoted by
\begin{equation}
u_{t_{0},x}:\R_+\to\X\quad,\quad
u_{t_{0},x}(t) =
\begin{cases}
X(t,t_{0})x &,\ t\geq t_{0}\\
0 &,\  0\leq t\leq t_0
\end{cases}
\quad.
\end{equation}

\begin{definition}
\label{def 2.2.}
An evolution family \Xts is said to be
\begin{enumerate}
\item[(u.e.s.)] \defnemph{uniformly exponentially stable}, if there exist $N, \nu > 0$ such that
\begin{gather}
\label{defn:expstab}
\|X(t,s)\|_{\lip} \leq Ne^{-\nu(t-s)}\ ,\ \mbox{for all}\ (t,s)\in\Delta\ ;
\end{gather}
\item[(u.s.)] \defnemph{uniformly stable} if there exists $N>0$ such that
\begin{equation}
\label{defn:unifstab}
\|X(t,s)\|_{\lip} \leq N\ ,\ \mbox{for all}\ (t,s)\in\Delta\ ;
\end{equation}
\item[(a.s.)] \defnemph{asymptotically stable} if all its trajectories are decaying to zero, i.e.
\begin{equation}
\label{defn:asymptstab}
\lim_{t\to\infty}{u_{t_0,x_0}(t)}=0 \ ,\ \mbox{for all}\ t_0\in\R_+\ \mbox{and}\ x_0\in\X\ .
\end{equation}
\end{enumerate}
\end{definition}

We will need in the next the following additional lemma.

\begin{lemma}
Let \Xts be a continuous evolution family on $\X$.
Then, the function $s\mapsto X(t,s)f(s):\R_+\to\X$ is locally integrable
provided that $f:\R_+\to\X$ is locally integrable.
\end{lemma}
\begin{proof}
This claim follows easily using conditions $(e_2)$, $(e_3)$ and $(e_4)$ of \Xts:
\begin{gather}
\begin{split}
\int_0^t \|X(t,s)f(s)-X(t,0)f(0)\|ds \leq &\, \int_0^t{\|X(t,s)f(s)-X(t,s)X(s,0)f(0)\|}ds \\
                \leq &\, Me^{\omega t}\int_0^t \left(\|f(s)\|+\|X(s,0)f(0)\|\right)ds
\end{split}
\end{gather}
which holds for any $t\in\R_+$.
\end{proof}

Given a continuous evolution family \Xts\!, we will denote by $\G$ the \emph{Green's operator}
\begin{equation}
\G:L^1_{loc}(\R_+,\X)\to L^1_{loc}(\R_+,\X)\quad,\quad (\G f)(t)=\int_0^t X(t,s)f(s)ds\quad.
\end{equation}
(where $L^1_{loc}(\R_+,\X)$ denotes the space of all locally integrable functions from $\R_+$ into $\X$).
Set now
$${\cal A}_{X} = \{ \chi _{[a,b]}u_{t_{0},x}:x\in \X,\; t_{0}\geq0,\; 0 \leq a \leq b \}\ ,$$
where $\chi_{[a,b]}$ denotes the characteristic function of the interval $[a,b]$.
Then, $A_X$ is contained in $L^p(\R_+,\X)$, for every $p\in[1,\infty]$.

\begin{definition}
\label{Definition 2.3.}
The pair $(L^{p}(\R_+,\X),L^{q}(\R_+,\X))$ is said to be \defnemph{admissible to \Xts}
if the following conditions hold
\begin{itemize}
        \item \quad the map $\G f$ lies in $L^q(\R_+,\X)$ for all $f\in L^{p}(\R_+,\X)$, and
        \item \quad $\G:L^{p}(\R_+,\X)\to L^{q}(\R_+,\X)$ is Lipschitz continuous, i.e. there exists $K>0$ such that
$$\|\G f-\G g\|_{q} \leq K\|f-g\|_{p}\ ,\ \mbox{for all}\ f,g\in L^p(\R_+,\X)\ .$$
\end{itemize}
\end{definition}

\section{Main results}

\begin{lemma}
\label{lemma:unifboundforh}
Let $h\in L^{q}(\R_+,\R)$, $q\in[1,\infty]$  such that $h(0)\geq0$. If
$$h(r)\leq m\,h(t)\quad,\ \mbox{for all}\ r\in[t,t+1]\ \mbox{and for all}\ t\geq0\ ,$$
then, $h\in L^{\infty}(\R_+,\R)$ and $\|h\|_{\infty}\leq mh(0)+ m\|h\|_{q}$.
\end{lemma}
\begin{proof}
Let $r\geq1$. Since $h(r)\leq m h(t)$ for all $t\in[r-1,r]$, it follows that
$$h(r)\,\leq\, m\int_{r-1}^r\!h(t)\,dt\,\leq\, m\|h\|_{q}\ .$$
If $r\in[0,1]$, from the hypothesis we have that $h(r)\leq m\,h(0)$.
Therefore,
$$h(r)\leq m (h(0) + \|h\|_{q})\ ,$$
for every $r\geq0$, and the above claim follows immediately.
\end{proof}

\begin{lemma}
\label{lem:masseraschaffer}
Let $g:\Delta\to\R_+$ be a function with the following properties:
\begin{enumerate}
\item $ g(t,t_{0}) \leq g(t,s)g(s,t_{0})$, for all $t \geq s \geq t_{0}$;
\item there exist $M, d > 0$ and $c \in (0,1)$ satisfying
\begin{gather*}
\begin{split}
& g(t,t_{0}) \leq M, \; \mbox{for all}\ t\in[t_{0},t_{0}+d]\,,\ t_0\geq0\ \mbox{and}\\
& g(t_{0}+d, t_{0})\leq c\,,\ \mbox{for all}\ t_{0}\geq 0\ .
\end{split}
\end{gather*}
\end{enumerate}
Then, there exist $N, \nu > 0$ such that
$$g(t,t_{0})\leq Ne^{-\nu (t-t_{0})}\ ,\ \mbox{for all}\; t \geq t_{0} \geq 0\ .$$
\end{lemma}
\begin{proof} Let $(t,t_0)\in\Delta$ and $n = \left[\frac{t-t_{0}}{d}\right]$.
Then, we have that
\begin{gather*}
\begin{split}
g(t,t_{0})\leq&\, g(t,t_{0}+nd)g(t_{0}+nd,t_{0})\leq
g(t,t_{0}+nd)c^{n} \\
\leq&\, Mc^{n}=Me^{-\nu nd}\leq Ne^{-\nu (t-t_{0})}\ ,
\end{split}
\end{gather*}
where $\nu = - \frac{1}{d}\ln c$, $N = Me^{\nu d}$.
\end{proof}

\medskip

Next we define the functions $a_{p}, b_{p}:\R_{+}\rightarrow \R$ given by
\[
a_{p}(t) = \left \{
\begin{array}{lll}
t^{1 - \frac{1}{p}} & , & p \in [1,\infty)\\
t & , & p = \infty
\end{array}
\right. ,
 b_{p}(t)= \| \chi _{[0,t]}\| _{p} = \left \{
\begin{array}{lll}
t^{\frac{1}{p}} & , & p \in [1,\infty) \\
1 & , & p = \infty .
\end{array}
\right.
\]

\begin{remark}
\label{rem 3.1.}
One can easily check that
$$\int\limits^{t_{0}+t}_{t_{0}} \| f(s)\| ds \leq a_{p}(t)\| f\| _{p}\ ,$$
for all $t,t_{0}\geq 0$ and $f \in L^{p}(\R_+,\X)$.
\end{remark}

\begin{theorem}
\label{the 3.2}
Let \Xts be a continuous evolution family and
$p,q\in[1,\infty]$ such that $(p,q)\neq(1,\infty)$.
If the pair $(L^{p}(\R_+,\X), L^{q}(\R_+,\X))$ is admissible to \Xts, then
\Xts is uniformly exponentially stable.
\end{theorem}
\begin{proof}
Let $t_{0}\geq 0,\; x_{1},\; x_{2} \in\X$ and
$f_{1},\; f_{2}: {{\R } }_{+}\rightarrow \X$ given by
\begin{gather}
f_{i}(t) =
\begin{cases}
X(t,t_{0})x_{i} &,\ t\in [t_{0}, t_{0}+1] \\
0 &,\ t\in\R_+\setminus [t_{0}, t_{0}+1]
\end{cases}
\quad,\ i ={1,2} \ .
\end{gather}
Clearly, $f_{1}, f_{2} \in {\cal A}_{X}$ with $\| f_{1}-f_{2}\|_{p}\leq Me^{\omega} \| x_{1}-x_{2}\| $
and
$$ (\G f_{i})(t) = \int\limits^{t}_{0}X(t,s)f_{i}(s)ds=\int\limits^{t_{0}+1}_{t_{0}}X(t,s)X(s,t_{0})x_{i}ds=
X(t,t_{0})x_{i}=u_{t_{0},x_{i}}(t)\ ,$$
for all $t\geq t_{0}+1$, $i ={1,2}$. For $t \in [t_{0}, t_{0}+1]$, we have that
$$\| u_{t_{0},x_{1}}(t)-u_{t_{0},x_{2}}(t)\| \leq
   \| X(t,t_{0})\|_{lip}\| x_{1}-x_{2}\| \leq Me^{\omega}\| x_{1}-x_{2}\|\ .$$
It follows that $u_{t_{0},x_{1}}-u_{t_{0},x_{2}} \in L^{q}(\R_+,\X)$ and
\begin{gather}
\begin{split}
\|u_{t_{0},x_{1}}-u_{t_{0},x_{2}}\|_{q} \leq&\, Me^{\omega} \|x_{1}-x_{2}\| + \|(\G f_{1}-\G f_{2})\,\chi _{[t_{0}+1,\infty)}\|_{q}\\
\leq&\, Me^{\omega}\|x_{1}-x_{2}\| + \| \G f_{1}- \G f_{2}\|_{q}\\
\leq&\, Me^{\omega}\|x_{1}-x_{2}\| + K \| f_{1}-f_{2}\|_{p}\\
\leq&\, (K+1)Me^{w}\| x_{1}-x_{2}\| .
\end{split}
\end{gather}
Let us define the map $h:\R_{+}\to\R_{+}$, $h(t) =\| u_{t_{0},x_{1}}(t_{0}+t) - u_{t_{0},x_{2}}(t_{0}+t)\|$.
Then, $h\in L^{q}(\R_+,\R)$ with
$\| h\|_{q} = \| u_{t_{0},x_{1}} -u_{t_{0},x_{2}}\|_{q} \leq (K+1)Me^{w}\| x_{1}-x_{2}\|$, and
\begin{gather*}
\begin{split}
h(r) =&\, \| X(t_{0}+r,t_{0}+t)X(t_{0}+t,t_{0})x_{1} - X(t_{0}+r,t_{0}+t)X(t_{0}+t,t_{0})x_{2}\| \\
   \leq&\, \| X(t_{0}+r,t_{0}+t)\| _{lip}\|X(t_{0}+t,t_{0})x_{1}-X(t_{0}+t,t_{0})x_{2}\| \\
   \leq&\, Me^{\omega(r-t)}h(t)\leq Me^{\omega}h(t)\quad,\quad 0\leq t\leq r\leq t+1\ .
\end{split}
\end{gather*}
By Lemma \ref{lemma:unifboundforh} we obtain that $h\in L^{\infty}(\R_+,\R)$ and
\begin{gather*}
\|h \| _{\infty} \leq Me^{\omega}\|h\| _{q} + Me^{\omega}h(0)\leq
 (K+1)M^{2}e^{2\omega}\|x_{1}-x_{2}\| + Me^{\omega}\|x_{1}-x_{2}\|
\end{gather*}
Now we have that there exists $C > 0$ such that
\begin{gather*}
\|X(t,s)x_{1}-X(t,s)x_{2}\| \leq C\| x_{1}-x_{2}\| \ ,
\end{gather*}
for all $(t,s)\in\Delta$ and  $x_{1}x_{2}\in \X$ and hence
\begin{gather}
\label{rel:unifstability}
\| X(t,s)\|_{\lip}\leq C\ , \quad \mbox{for all}\quad (t,s)\in\Delta\ .
\end{gather}

Consider again $x_{1},x_{2}\in \X$, $t_{0}\geq 0$, $\delta > 0$, $f_{3}, f_{4}:\R_{+}\rightarrow \X$
given by
\begin{gather}
f_{i}(t)=
\begin{cases}
X(t,t_{0})x_{i-2} &,\ t\in [t_{0}, t_{0}+\delta]  \\
0 &,\  t\in\R_{+}\setminus [t_{0},t_{0}+\delta ]
\end{cases}
\quad,\ i = {3,4}\ .
\end{gather}
Then $f_{3}, f_{4}\in {\cal A}_{X}$ with $\| f_{3}-f_{4}\| _{p}\leq C b_{p}(\delta)\| x_{1}-x_{2}\|$
and
$$
(\G f_{i})(t) = \int\limits^{t}_{0}X(t,s)f_{i}(s)ds=
        \int\limits^{t}_{t_{0}}X(t,s)X(s,t_{0})x_{i-2}\,ds =
        (t-t_{0})U(t,t_{0})x_{i-2}
$$
for all $t\in[t_{0},t_{0}+\delta]$, $i ={3,4}$.
Using the last equality we have that
\begin{gather}
\begin{split}
\frac{\delta ^{2}}{2}\| X(t_{0}+\delta,t_{0})x_{1} -X(t_{0}+\delta,t_{0})x_{2}\| =&\,
        \int\limits^{t_{0}+\delta}_{t_{0}}(t-t_{0})\|X(t_{0}+\delta,t_{0})x_{1}-X(t_{0}+\delta, t_{0})x_{2}\|dt\\
   \leq&\, C\int\limits^{t_{0}+\delta}_{t_{0}}(t-t_0)\|X(t,t_{0})x_{1}-X(t,t_{0})x_{2}\| dt\\
   =&\,C\int\limits^{t_{0}+\delta}_{t_{0}}\| (\G f_{3})(t)- (\G f_{4})(t)\| dt\\
  \leq&\, C a_{q}(\delta) \| \G f_{3}- \G f_{4}\|_{q}\\
  \leq&\, KCa_{q}(\delta)\|f_{3}-f_{4}\|_{p}\\
  \leq&\, KC^{2}a_{q}(\delta)b_{p}(\delta)\|x_{1}-x_{2}\|\\
        =&\, \frac{KC^{2}\delta^{2}}{a_{p}(\delta)b_{q}(\delta)}\|x_{1}-x_{2}\|\ .
\end{split}
\end{gather}
It follows that
$$\|X(t_{0}+\delta, t_{0})\| _{lip}\leq\frac{2KC^{2}}{a_{p}(\delta)b_{q}(\delta)}\ ,$$
for all $t_{0}\geq 0$, $\delta > 0$.
Since  $(p,q)\not=(1,\infty)$, we have that $\lim\limits_{\delta \to\infty}a_{p}(\delta)b_{q}(\delta)=\infty$,
and therefore we can choose $d>0$ such that
\begin{gather}
\|X(t_{0}+d,t_{0})\|_{\lip}\leq \frac{1}{2}\ ,
\end{gather}
{for all} $t_{0}\geq0$.
Applying Lemma~\ref{lem:masseraschaffer} to the map
$g:\Delta\rightarrow \R_{+}$, defined by $ g(t,s)=\| U(t,s)\|_{lip}$
we obtain that \Xts is uniformly exponentially stable.

\end{proof}

\begin{theorem}\label{the 3.1}
%Let \Xts be a continuous evolution family.
%
%Then,
The pair $(L^{1}(\R_+,\X),L^{\infty}(\R_+,\X))$
is admissible to the continuous evolution fa\-mi\-ly \Xts
if and only if the following statements hold
\begin{enumerate}
\item[$(i)$]\quad there exists $\psi\in L^1(\R_+,\X)$ such that $\G\psi\in L^{\infty}(\R_+,\X)$, and
\item[$(ii)$]\quad there exists $N>0$ such that $\|X(t,s)\|_{lip} \leq N$, for all $(t,s)\in\Delta$.
\end{enumerate}
\end{theorem}
\begin{proof}
The {\it necessity} follows from the proof of Theorem~\ref{the 3.2} since we
proved \eqref{rel:unifstability} without the $(p,q)\neq(1,\infty)$ assumption.

{\it Sufficiency}. Let $f\in L^{\infty}(\R_+,\X)$.
Since $\psi\in L^1(\R_+,\X)$ and $f\in L^1(\R_+,\X)$, we have that
\begin{gather}
\begin{split}
\label{unifstabil:estim1}
\| (\G\psi)(t)-(\G f)(t)\| \leq&\, \int\limits^{t}_{0}\| X(t,s)\phi(s)- X(t,s)f (s)\| ds  \\
\leq&\, \int\limits^{t}_{0}\|X(t,s)\|_{lip}\|\psi(s) - f(s)\|ds \\
\leq&\, N\int\limits^{t}_{0}\| \psi (s) -f (s)\| ds \\
\leq& N \| \psi  - f  \| _{1} \quad,
\end{split}
\end{gather}
for all $t \geq0$, which implies that $\G f \in L^\infty (\R_+,\X)$.

It remains to prove that $\G: L^1(\R_+,\X)\to L^{\infty}(\R_+,\X)$ is Lipschitz continuous.
But, using similar arguments as in \eqref{unifstabil:estim1}, we obtain that
\begin{gather}
\|\G f_1- \G g\|_{L^{\infty}(\R_+,\X)}\leq N\|f-g\|_{L^1(\R_+,\X)}\ ,
\end{gather}
for all $f,g \in L^1(\R_+,\X)$.
\end{proof}

\begin{remark}
Theorem \ref{the 3.1} extends a similar result for linear
equations (see e.g. \cite{cop,masseraschaffer}). Note that in the linear
case, condition $(ii)$ is automatically satisfied.
\end{remark}

In the following theorem we try to answer concerns
regarding the converse of what was obtained in Theorem~\ref{the 3.2}.
With elementary arguments we can show that the admissibility of the pair
$(L^p(\R_+,\X),L^q(\R_+,\X))$ with $p\leq q$ is a necessary condition for
the uniform exponential stability of a continuous evolution family.
The idea is based on the use of Fubini's theorem, H\"{o}lder's inequality
and the observation that if $f\in L^p(\R_+,\X)\cap L^q(\R_+,\X)$,
then $f\in L^r(\R_+,\X)$ with $\|f\|_r\leq\max\{\|f\|_p,\|f\|_q\}$,
for any $1\leq p\leq r\leq q\leq\infty$.

\begin{theorem}
Let \Xts be a continuous evolution family and $1\leq p\leq q\leq \infty$.

Then, the pair $(L^p(\R_+,\X), L^q(\R_+,\X))$ is admissible to \Xts provided that
\begin{enumerate}
\item there is a function $\psi\in L^p(\R_+,\X)$ such that $\G\psi \in L^q(\R_+,\X)$, and
\item \Xts is uniformly exponentially stable.
\end{enumerate}
\end{theorem}
\begin{proof}
Let $N,\nu>0$ be like in Definition~\ref{def 2.2.}.
Fix arbitrarily $f,g\in L^p(\R_+,\X)$. We have that
\begin{eqnarray*}
\| (\G f)(t) -(\G g)(t)\| &\leq& \int_0^t \|X(t,s)(f(s)-g(s))\|ds \leq \int_0^t \|X(t,s)\|_{lip}\|f(s)-g(s)\|ds\\
 &\leq& N\int_0^t e^{-\nu(t-s)}\|f(s)-g(s)\|ds
\end{eqnarray*}
Consider $h,H:\R_+\to\R_+$ given by $h(t)=\|f(t)-g(t)\|$ and
\begin{equation}
H(t)=\int_0^t e^{-\nu(t-s)}h(s)ds\
\end{equation}
for any $t\in\R_+$.
In what follows, we will prove that if $h\in L^p(\R_+,\R)$ then $H\in L^q(\R_+,\R)$
(in the hypothesis $1\leq p\leq q\leq\infty$).

\textit{Case 1.} If $p=\infty$, then $q=\infty$ and since
\begin{eqnarray*}
H(t) &\leq& \int_0^t e^{-\nu(t-s)}\|h\|_{\infty}ds\leq \|h\|_{\infty}\int_0^t e^{-\nu\tau}d\tau \leq \frac{1}{\nu}\|h\|_{\infty}\ ,
\end{eqnarray*}
for all $t\in\R_+$, it follows that $H\in L^{\infty}(\R_+,\R)$ with $\|H\|_{\infty}\leq \frac{1}{\nu}\|h\|_{\infty}$.

\textit{Case 2.} If $p=1$, note that
\begin{eqnarray*}
H(t) = \int_0^t e^{-\nu(t-s)}h(s)ds\leq \int_0^t h(s)ds \leq \|h\|_1\ ,
\end{eqnarray*}
for all $t\in\R^+$ and thus $H\in L^{\infty}(\R_+,\R)$ with $\|H\|_{\infty}\leq \|h\|_1$. Also, using Fubini's theorem we have that
\begin{eqnarray*}
\int_0^{\infty} H(t)dt &=& \int_0^{\infty}\!\int_0^{t}e^{-\nu(t-s)}h(s)\,ds\,dt =
  \int_0^{\infty}\!\int_s^{\infty} e^{-\nu(t-s)}h(s)\,dt\,ds \\
  &=& \int_{0}^{\infty}e^{\nu s}h(s)\int_s^{\infty}\!e^{-\nu t}dt\,ds= \int_0^{\infty}e^{\nu s}h(s)\frac{e^{-\nu s}}{\nu}ds\\
  &=& \frac{1}{\nu}\|h\|_1\ ,
\end{eqnarray*}
which implies that $H\in L^1(\R_+,\R)$ with $\|H\|_1\leq \frac{1}{\nu}\|h\|_1$.
Then, $H\in L^q(\R_+,\R)$ and $\|H\|_q\leq \max\{1,1/\nu\} \|h\|_1$.

\textit{Case 3.} If $p\in(1,\infty)$, let $p'\in(1,\infty)$ such that $\frac{1}{p}+\frac{1}{p'}=1$ and
let $\alpha,\beta\in(0,1)$ such that $\alpha+\beta=1$. For any $t\in\R_+$, we can write down
\begin{eqnarray*}
H(t)\leq \left(\int_0^t h(s)^p ds\right)^{1/p} \left(\int_0^t e^{-\nu p'\tau}d\tau\right)^{1/p'}\leq \frac{1}{(\nu p')^{1/p'}}\|h\|_p\ .
\end{eqnarray*}
We obtained that $H\in L^{\infty}(\R_+,\R)$ and $\|H\|_{\infty}\leq (\nu p')^{-1/p'}\|h\|_p$.
Next, we prove that $H\in L^p(\R_+,\R)$.
Via H\"{o}lder's inequality, we have that
\begin{eqnarray*}
\int_0^t e^{-\nu(t-s)}h(s)ds &=& \int_0^t e^{-\nu\alpha(t-s)}e^{-\nu\beta(t-s)}h(s)ds \\
  &\leq& \left(\int_0^t e^{-\nu\alpha p'(t-s)}ds\right)^{1/p'}\left(\int_0^t e^{-\nu\beta p(t-s)}h(s)ds\right)^{1/p}\\
  &\leq& \left[\left(\frac{1}{\nu\alpha p'}\right)^{p-1} \int_0^t e^{-\nu\beta p(t-s)}h(s)ds\right]^{1/p}\ .
\end{eqnarray*}
Then, denoting $C:=(\nu\alpha p')^{1-p}>0$, we can write down
\begin{eqnarray*}
\int_0^{\infty}H(t)^pdt &=& \int_0^{\infty}\left(\int_0^t e^{-\nu(t-s)}h(s)ds\right)^p dt \leq C
  \int_{0}^{\infty}\int_0^{t} e^{-\nu\beta p(t-s)}h(s)^pds\,dt \\
  &\leq& \frac{C}{\nu\beta p} \|h\|_p^p
\end{eqnarray*}
(the last step follows similarly to \textit{Case 2}, using Fubini's theorem).
From here, it follows that $H\in L^p(\R_+,\R)$ with $\|H\|_p\leq C^{1/p}(\nu\beta p)^{-1/p}\|h\|_p$.
Therefore,
$H\in L^q(\R_+,\R)$ and we have that
$\|H\|_q\leq\max\{(\nu p')^{-1/p'},C^{1/p}(\nu\beta p)^{-1/p}\}\|h\|_p$.

In any case, we obtained that $H\in L^q(\R_+,\R)$ (provided that $h\in L^p(\R_+,\R)$) and
the existence of some $K>0$ (independent of $h$) such that $\|H\|_q\leq K\|h\|_p$.

To complete the proof, note that from (a) we have that $\psi-f\in L^p(\R_+,\X)$ and $\G\psi\in L^q(\R_+,\X)$;
in virtue of all above we obtain that $\G\psi - \G f\in L^q(\R_+,\X)$
(no matter how we take $f\in L^p(\R_+,\X)$).
Moreover, we can state that
$\|\G f - \G g\|_q\leq K\|f-g\|_p$, for all $f,g\in L^p(\R_+,\X)$.
Hence, $\G\in\Lip(L^p(\R_+,\X),L^q(\R_+,\X))$.

\end{proof}

\begin{theorem}
If the pair $(L^1(\R_+,\X),L_0^{\infty}(\R_+,\X))$ is admissible
to  the continuous evolution family \Xts\!, then \Xts is asymptotically stable.
\end{theorem}
\begin{proof}
Let $x\in\X$, $t_0\in\R_+$ and consider the function $f:\R_+\to\X$
$$f(t) =
\begin{cases}
X(t,t_0)x &,\ t\in[t_0,t_0+1]\\
0 &,\ t\in\R_+\setminus[t_0,t_0+1]
\end{cases}\quad.$$
We have that $f\in L^1(\R_+,\X)$ and note that
$$(\G f)(t) =\int_{t_0}^{t_0+1}X(t,s)X(s,t_0)xds=u_{t_0,x}(t)\ ,$$
for all $t\geq t_0+1$. Since $\G f\in L_0^{\infty}(\R_+,\X)$, we get that $\lim_{t\to\infty}u_{t_0,x}(t)=0$.

\end{proof}

\section{Applications}
In this section as a model for applications of the obtained
results in the previous section we consider equations of the form
\begin{eqnarray*}
\frac{\partial u(t,x)}{\partial t} &=&  \frac{\partial^2 u(t,x)}{\partial x^2} +g(t,u(t,x)), \quad t>0 ,x\in (0,\pi ),\\
\frac{\partial u(t,x)}{\partial t} &=& 0, \quad x=0,\pi ,
\end{eqnarray*}
where $u(t,x)$ is a scalar function of $(t,x)\in \R^+\times\R$,
$g(t,y)$ is uniformly Lipschitz continuous in $y\in\R$, and
$g(\cdot , 0)\in L^r(\R_+,\X)$ with $1\le r < \infty$. This equation can be
re-written in the following abstract form in a Banach space $ \X$
\begin{equation}\label{ex1}
\frac{d}{dt} u(t)=  \Delta u(t) + G(t,u(t)),
\end{equation}
where $X= \{ v\in W^{2,2}(0,\pi ): v' =0 \ \mbox{at} \ x=0,\pi
\}$, $ \Delta v=  \partial^2/\partial x^2$ on $X$,  $G(t,u(t))=
g(t,u(t,\cdot ))$. As is well known, (actually, an extension of)
$\Delta$ generates a strongly continuous analytic semigroup in $\X$
that we denote by $(T(t))_{t\ge 0}$. By a standard argument (see
e.g. \cite{aulmin}), we can prove that (\ref{ex1}) generates a
continuous evolution family $\Xts$ in $\X$ that is determined from the equation
\begin{equation}
X(t,s)x= T(t-s)x + \int^t_s T(t-\xi ) G(\xi ,X(\xi ,s)x)d\xi ,
\quad \ \mbox{for all} \ t\ge s \ge 0.
\end{equation}
\begin{theorem}\label{the 3.3}
Assume that the pair $(L^p(\R_+,\X),L^q(\R_+,\X))$ is admissible to the continuous evolution family $\Xts$.
Then there exists a mild solution $u\in L^r(\R_+,\X)$
of (\ref{ex1}) that attracts all other mild solutions of the
equations at exponential rate.
\end{theorem}
\begin{proof}
Let us consider the evolution semigroup $(T^h)_{h\ge 0}$ in $L^r(\R_+,\X)$
associated with the linear equations $u'=\Delta u$, defined as
\begin{equation}\label{semi}
[T^hf](t):= \begin{cases} T(h)f(t-h), \ \mbox{if} \ h\le t \\
0, \ \mbox{if} \ 0\le t< h ,
\end{cases}
\end{equation}
for all $f\in L^r(\R_+,\X)$. Since $1\le r<\infty$, this semigroup is
strongly continuous (see e.g. \cite{chilat,minrabsch}). Let us
denote by ${\cal L}$ the generator of this semigroup. As is well
known (see e.g. \cite{minrabsch}), $u\in D({\cal L})$ and ${\cal
L}u=-f$ if and only if
\begin{equation}\label{ex 2}
u(t)= \int^t_0 T(t-\xi )f(\xi d\xi \quad t\ge 0.
\end{equation}
Consider the operator ${\cal L} +{\cal G}$ on $L^r(\R_+,\X)$, where ${\cal
G}$ is the Nemytsky operator associated with $G$, defined as $L^r(\R_+,\X):
\ni \phi (\cdot ) \mapsto G(\cdot , \phi (\cdot )) \in L^r(\R_+,\X)$. Note
that under the assumption on $g$, the Nemytsky operator acts in
$L^r(\R_+,\X)$ as a Lipschitz continuous operator. So, in the same way as
in \cite{aulmin}, we can show that the operator ${\cal L} +{\cal
G}$ associated with (\ref{ex1}) generates a strongly continuous
semigroup $(S(h))_{h\ge 0}$ in $L^r(\R_+,\X)$ that is referred to as the
evolution semigroup associated with (\ref{ex1}). Moreover, this
semigroup is determined by
\begin{equation}\label{semi-1}
[S(h)f](t):= \begin{cases} X(t,t-h)f(t-h), \ \mbox{if} \ h\le t \\
0, \ \mbox{if} \ 0\le t< h .
\end{cases}
\end{equation}
for all $f\in L^r(\R_+,\X)$. By the admissibility of the pair $(L^p(\R_+,\X),L^q(\R_+,\X))$,
$\Xts$ is exponentially stable. This implies the strict
contraction of $S(h)$ for sufficiently large $h$. In fact,
\begin{eqnarray*}
\| S(h)\phi -S(h)\psi \| _r &\le& Ne^{-\alpha h} \| \phi -\psi\|_r
\end{eqnarray*}
for all $\phi, \psi \in L^r(\R_+,\X)$. Therefore, if $h$ is sufficiently
large $Ne^{-\alpha h} <1$. Let us fix a sufficiently large integer
$n_0$. This yields that $S(n_0)$ has a unique fixed point
$\varphi\in L^r$. Since $S(h)$ commutes with $S(n_0)$ for all
$h\ge 0$, it is easy to see that $\varphi$ is the unique common
fixed point for the entire semigroup $(S(h))_{h\ge 0}$. This
implies that $({\cal L}+{\cal G})\varphi =0$. So, ${\cal L}\varphi
= -{\cal G}\varphi$, and by the formula (\ref{ex 2}), we have
$$
\varphi (t)=\int^t_0 T(t-\xi )G(\xi ,\varphi (\xi ))d\xi ,\quad
t\ge 0 .
$$
This means that $\varphi\in L^r(\R_+,\X)$ is a mild solution starting at
zero of (\ref{ex1}), so $\varphi (t) =X(t,0)0.$ Now we show that
this solution attracts all other solutions at exponential rate. In
fact, every other solution starting, say, at $x\in \X$ is of the
form $X(t,0)x$. Therefore,
$$
\| X(t,0)x-X(t,0)0\| \le Ne^{-\alpha t} \| x\| , \quad t\ge 0 .
$$
This completes the proof of the theorem.
\end{proof}
Before concluding this section we give an application of Theorem
\ref{the 3.1}
\begin{proposition}
Let the pair $(L^1(\R_+,\X),L^\infty(\R_+,\X))$ be admissible to the continuous evolution family $\Xts$ (generated by (\ref{ex1})). Moreover, assume
that $g(t,0)=0$ for all $t\ge 0$. Then every mild solution $u$ of
(\ref{ex1}) is bounded.
\end{proposition}
\begin{proof}
Obviously, the $u(t)=0$ is the trivial mild solution of
(\ref{ex1}). Therefore, for every mild solution $u(t) = X(t,s)x$,
where $x\in \X$, we have
$$
\| u(t)\| = \| X(t,s)x\| =\| X(t,s)x-X(t,s)0\| \le \| X(t,s)\|
_{lip} \| x\| \le N \| x\| \quad \ \mbox{for all} \ t\ge 0,
$$
where $N$ is a positive number whose existence is guaranteed by
Theorem \ref{the 3.1}.
\end{proof}

\bibliographystyle{amsplain}

\end{document}